\documentclass[]{article}
\usepackage[hidelinks,draft=false]{hyperref}
\usepackage{times}
\usepackage{shuffle}
\usepackage{booktabs}
\usepackage{microtype}
\usepackage[final]{graphicx}
\usepackage{geometry}
\usepackage{comment}
\usepackage{mathtools} 
\usepackage{subfig}
\usepackage{float}
\usepackage{amsmath}
\makeatletter
\let\over\@@over
\let\atop\@@atop
\makeatother
\usepackage{amssymb}
\usepackage{mathrsfs}
\usepackage{mhequ}
\usepackage{./styles/mhenvs}
\usepackage{./styles/mhfig}
\usepackage{./styles/mhsymb}

\usepackage{hyperref}

\usepackage{times}
\usepackage{microtype}
\usepackage{amsmath}
\usepackage{amssymb}
\usepackage{tikz}
\usepackage{wasysym}
\usepackage{centernot}

\colorlet{symbols}{blue!90!black}

\def\1{\mathbf{1}}

\colorlet{testcolor}{green!60!black}

\usetikzlibrary{shapes.misc}
\usetikzlibrary{shapes.symbols}
\usetikzlibrary{decorations}
\usetikzlibrary{decorations.markings}

\makeatletter
\def\DeclareSymbol#1#2#3{\expandafter\gdef\csname MH@symb@#1\endcsname{\tikz[baseline=#2,scale=0.15]{#3}}%
\expandafter\gdef\csname MH@symb@#1s\endcsname{\scalebox{0.6}{\tikz[baseline=#2,scale=0.15]{#3}}}}
\def\<#1>{\csname MH@symb@#1\endcsname}
\makeatother

\DeclareSymbol{X}{-2.4}{\node[circle, fill, inner sep=1pt] {};}

\DeclareSymbol{1}{0}{\draw[white] (-.4,0) -- (.4,0); \draw (0,-.2)  -- (0,1.6) node[dot, fill, inner sep=0.7pt] {};}

\newop{diag}

\newcommand\minus{%
  \setbox0=\hbox{-}%
  \vcenter{%
    \hrule width\wd0 height \the\fontdimen8\textfont3%
  }%
}

\def\E{\mathbb{E}}
\def\P{\mathbb{P}}
\def\/{\,\rule[-0.25em]{2pt}{1em}\,}

\def\${|\!|\!|}

\def\dd{\mathrm{d}}

\def\R{\mathbb{R}}

\def\N{\mathbb{N}}

\newcommand{\norm}[1]{\left\lVert #1 \right\rVert}

\DeclareMathOperator*{\weaklim}{\mathrm{w}\text{-}\mathrm{lim}}

\usepackage{xcolor}
\usepackage{tcolorbox}

\begin{document}  

\title{A proof of the Lyons-Pemantle-Peres monotonicity conjecture for high offspring laws}

\author{Felix Schwarzfischer}

\maketitle 

\begin{abstract}
    We prove the Lyons-Pemantle-Peres conjecture for high enough offspring laws, that is we show that the speed satisfies $\ell'(\lambda) < 0 $ for the entire transient regime of the bias $\lambda  \in (0, m)$ in the special case that the number of children $Z$ in the Galton-Watson tree satisfies $Z \geq 4$ almost surely. 
\end{abstract}

\tableofcontents

\section{Introduction}
Let $Z$ be a $\N = \{1,2,\dots \}$ valued random variable and denote by $\CT$ the corresponding (leafless) Galton-Watson tree having $Z$ as the offspring law at each vertex $v$ and call the resulting law of $\CT$ as $\mathbf{P}$. On this random environment, we then introduce a random walk $ (\xi^\lambda_k)_{k \in \N_0}$ starting at the root $\rho $ of $\CT$, and which follows the following transition probabilities for each vertex $v$ of $\CT$ which is not the root $\rho$:
\begin{equs}
P_{\lambda, \CT}(\xi^\lambda_{k+1}= \mathrm{parent}(v) | \xi^\lambda_k  = v  ) &= \frac{\lambda}{\lambda + Z(v)} \\ 
 P_{\lambda, \CT}(\xi^\lambda_{k+1} = \mathrm{child}_i(v) | \xi^\lambda_k = v ) &= \frac{1}{\lambda + Z(v) } , \quad  i=1 ,\dots, Z(v) 
\end{equs}
as well as for $v = \rho$ 
\begin{equ}
    P_{\lambda, \CT}(\xi^\lambda_{k+1} = \mathrm{child}_i(\rho)| \xi^\lambda_k = \rho) = \frac{1}{Z(\rho)} \, . 
\end{equ}
Here $\lambda \in \R_{\geq 0} $ is called the bias parameter, indicating the bias to return to the parent nodes during the random walk. Going from the quenched law $P_{\lambda, \CT}$, we write for the annealed probability measure 
\begin{equ}
(\P \circ (\CT, \xi^\lambda)^{-1})(\dd \tau, \dd x)  = \mathbf{P}(\dd \tau) P_{\lambda, \tau}(\dd x) \, . 
\end{equ}
Here we model everything on one probability space $(\Omega, \CF, \P)$ furnishing $ \CT, \xi^\lambda$ and any additional random variables as needed on the fly. 
One can prove that for $\lambda \in (0,m)$ with $m := \E Z$, the asymptotic speed 
\begin{equ}
\ell (\lambda) := \lim_{k \to \infty} \frac{|\xi^\lambda_k|}{k} 
\end{equ}
exists and is deterministic $\P$-a.s. \cite[Theorem 3.1 + Section 4]{lyons1996biased}. 
Given this context, we are going to prove the following. 
\begin{theorem}\label{thm:lpp_four_children}
    Let $Z$ be such that $Z \geq 4$ almost surely and $m = \E Z < \infty $. Then $\ell: (0,m) \to \R_{\geq 0}$ is differentiable and 
    \begin{equ}
    \forall \lambda \in (0,m) : \ell'(\lambda) < 0 \, . 
    \end{equ}
\end{theorem}
In recent years variants of Theorem \ref{thm:lpp_four_children} have been treated, providing partial resolutions to the Lyons-Pemantle-Peres monotonicity conjecture \cite{Lyonsetal1997}. We mention for example \cite{arous2011proof} and \cite{SongWangXiang2025}, where $\ell'(\lambda) < 0$ was shown for $\lambda \leq \frac{m_1}{1160}$ and $\lambda \leq \frac{m_1}{1+ \sqrt{1-\frac{1}{m_1}}}$ respectively, with $m_1 = \min\{ k \geq 1: \P(Z = k) > 0 \} $. Our advance is the extension to the entire transient regime $\lambda \in (0,m)$, albeit at the cost of requiring high off spring laws $Z \geq 4$. With similar techniques as presented in this note, it is possible to prove Theorem \ref{thm:lpp_four_children} for the less restrictive case $Z \geq 2$.  In the interest of clarity we refrain from doing so here, as this requires just lengthier arguments that are merely technical in nature which nonetheless don't appear to be sufficient yet to resolve the $Z \geq 1$ case. 
    
To comment on the main technical concept appearing in the proof, recall the tree-dependent escape probability 
\begin{equ}
 X_\lambda(\CT  ) = \P(\xi^\lambda \textrm{ never returns to  }\rho^* | \CT ) = P_{\lambda, \CT}(\xi \textrm{ never returns to  } \rho^* ) \, . 
\end{equ}
Here, $\rho^*$ is an artificial parent of the root $\rho$, and the random walk $(\xi^\lambda_k)$ starting at $\rho$ on the extended tree $\CT_*$ is similarly defined as the ordinary random walk on $\CT$, resulting though in the same speed $\ell(\lambda)$ as in the original random walk \cite[Lemma 1]{SongWangXiang2025} for $\lambda \in (0,m)$. 
Calling the law of $X_\lambda$ under $\P$ as $\mu_\lambda$, we recall the nonlinear fixed point equation satisfied by the law $\mu_\lambda$ expressed via $X_\lambda$ and $Z$ as 
\begin{equ}\label{eq:fixed_point_law}
X_\lambda \stackrel{d}{=} \frac{\sum_{i=1}^Z X_{\lambda, i}}{\lambda + \sum_{i=1}^Z X_{\lambda, i}} \, , 
\end{equ}
where $X_{\lambda, i}$ are i.i.d. with law $\mu_\lambda$, and equality is in distribution under $\P$. This equation \eref{eq:fixed_point_law} is  explicitly mentioned for example in \cite[Section 13.3.3]{aldous2002reversible}, but can essentially be traced back to \cite{lyons1995ergodic,lyons1996biased}. Whereas \cite{SongWangXiang2025} now differentiates in a pathwise sense the the random variable $X_\lambda$ in \eref{eq:fixed_point_law}\footnote{not literally \eref{eq:fixed_point_law} but an analogous equation that holds pathwise (cf. \cite[Section 2]{SongWangXiang2025})}, we differentiate the \emph{law} $\mu_\lambda$, and thus derive the linearized equation around a given $\lambda$ for a \emph{measure} $\Gamma_\lambda$, resulting in an instance of the Dyson equation 
\begin{equ}\label{eq:fixed_point_for_Gamma_intro}
\Gamma_\lambda = \sigma_\lambda  + T_\lambda  \Gamma_\lambda \, . 
\end{equ}
Let us elaborate more in an informal manner how to arrive at \eref{eq:fixed_point_for_Gamma_intro}. For every test function $f$, we can write \eref{eq:fixed_point_law} as an explicit fixed point equation for $\mu_\lambda$: 
\begin{equ}\label{eq:fixed_point_law_explicit_against_f}
\int f(x) \dd \mu_\lambda(x) = \sum_n p_n \int f\left(\frac{\norm{x}_1}{\lambda + \norm{x}_1 }\right) \dd \mu_\lambda^{\otimes n}(x) 
\end{equ}
where $p_n= \P(Z = n)$ and $\norm{x}_1 = x_1+ \dots + x_n$. We want to study now the first order change to $\mu_\lambda$ when  $\lambda \mapsto \lambda + \eps $, that is, we linearize $\mu_{\lambda+\eps} \approx \mu_\lambda + \eps \dot{\mu_\lambda}$. Plugging into \eref{eq:fixed_point_law_explicit_against_f}, to order $\eps$, we get easily the equation satisfied by $\dot{\mu}_\lambda$: 
\begin{equ}\label{eq:fixed_point_mu_dot}
\int f(x) \dd \dot{\mu}_\lambda = - \frac{1}{\lambda} \E f'(X_\lambda)X_\lambda(1-X_\lambda) + \sum_n n p_n \int f\left( \frac{\norm{x}_1}{\lambda + \norm{x}_1} \right) \dd \left( \mu_\lambda^{n-1} \otimes \dot{\mu}_\lambda  \right)(x) \, . 
\end{equ}
On the other hand, for the family $(\mu_\lambda)_{\lambda \in (0,m)}$, pretending $\lambda $ to be a time parameter, we have a fairly general `continuity equation' of the form 
\begin{equ}
\dot{\mu}_\lambda(\dd x) = \partial_\lambda \mu_\lambda(\dd x) = \partial_\lambda \left[ \frac{F_\lambda(x+ \dd x) - F_\lambda(x)}{\dd x} \dd x \right]  = \partial_x \partial_\lambda F_\lambda(x) \dd x =: \partial_x \Gamma_\lambda(\dd x)  \, . 
\end{equ}
where $F_\lambda(x) = \mu((-\infty, x])$. This terminology is chosen, as it resembles the continuity equation from fluid dynamics $\partial_t \rho  + \nabla v \rho = 0 $, with $\Gamma_\lambda$ representing the `velocity field' in this context. This is to say, that it is more natural to formulate \eref{eq:fixed_point_mu_dot} in terms of $\Gamma_\lambda$ than with $\dot{\mu}$. Hence, with formal integration by parts
\begin{equ}
\int f(x) \dot{\mu}(\dd x) = \int f(x)  \partial_x \Gamma_\lambda(\dd x) = - \int f'(x) \Gamma_\lambda(\dd x) \, , 
\end{equ} 
we arrive at, by calling the measure arising in the first term of \eref{eq:fixed_point_mu_dot} as $\sigma_\lambda$, 
\begin{equ}
\int f'(x) \Gamma_\lambda(\dd x) = \int f'(x) \sigma_\lambda(\dd x) + \sum_n n p_n \int f'\left( \frac{\norm{x}_1}{\lambda + \norm{x}_1}\right) \frac{\lambda}{(\lambda + \norm{x}_1)^2} \dd(\mu_\lambda^{n-1} \otimes \Gamma_\lambda)(x)  \, . 
\end{equ}
As this holds for every $f$ nice enough, we can read off the operator $T_\lambda$ and arrive at \eref{eq:fixed_point_for_Gamma_intro}. 

This entire preprocessing step lets us eventually express the derivative $\ell'(\lambda)$ in a way where its sign can be easily read off from evaluations of some bilinear form on explicit measures related to $\mu_\lambda$ and $\Gamma_\lambda$. 
\paragraph{Notation:} We will frequently use $X_\lambda(\CT)$ as the escape probability of the random walk on a given tree $\CT$ starting at the root. We will drop $\CT$ though from the notation. Furthermore, many of the appearing objects $X_\lambda, \mu_\lambda, Q_\lambda, T_\lambda$ will be dependent on $\lambda$, however for legibility's sake, inside the proofs we will often drop the subscript $\lambda$. 
\section{A formula for the speed $\ell'(\lambda)$ }
We first recall from Aidekon \cite{aidekon2014speed} and Rousselin \cite{Rousselin2018} the formula for the speed 
\begin{equ}\label{eq:rousselin_formula}
\ell(\lambda) = \frac{A_\lambda}{B_\lambda} \, , 
\end{equ}
with  
\begin{equ}
A_\lambda := \E \frac{XY}{D_\lambda(X,Y)} \, , B_\lambda := \E \frac{X+Y - XY}{D_\lambda(X,Y)}  \, , 
\end{equ}
where $D_\lambda(x,y) = \lambda - (1-x)(1-y)$. Here $X,Y$ are i.i.d and $X \sim \mu_\lambda = \mathrm{Law}_{\P}(X_\lambda)$. 
It turns out that for the technical details appearing later on, it is very convenient to have knowledge of where $\mu_\lambda$ is supported, since this gives control on the kernel $D_\lambda$. 
\begin{lemma}
    For any offspring law such that $Z \geq 4$, it holds that 
    \begin{equ}
    X_\lambda \geq \left(1-\frac{\lambda}{4}\right)_+ \, ,  \P\mathrm{-a.s.} 
    \end{equ}
    and hence $\supp \mu_\lambda \subset \left[\left(1-\frac{\lambda}{4}\right)_+, 1 \right] $. As a consequence, we get two-sided bound on the kernel $D_\lambda$ when sampling from the distribution $\mu_\lambda$:
    
    \begin{equ}\label{eq:D_kernel_bounds}
     \forall x,y \in \supp\mu_\lambda :  \,  \frac{3}{4} \lambda \leq D_\lambda(x,y) \leq \lambda 
    \end{equ}

\end{lemma}
\begin{proof}
    This is a special case of \cite[Lemma 3]{SongWangXiang2025}, applied to $m_1 =: \inf\{ n \geq 1 : \P(Z =n ) > 0 \}  \geq  4$, with \eref{eq:D_kernel_bounds} an immediate computation. 
\end{proof}
It will be useful to introduce short-hand notations for the following sets of measures 
\begin{equs}
\CM &:= \{ \nu | \nu  \textrm{ is signed measure on } ([0,1], \CB([0,1])) \}   \\ 
\CM_\lambda &:= \{ \nu \in \CM | \supp \nu \subset [(1-\lambda / 4)_+ , 1] \} \\ 
\CM^+ &:= \{ \nu \in \CM | \nu \geq 0 \} \, \\ 
\CM^+_\lambda &:= \CM_\lambda \cap \CM^+ \, . 
\end{equs}
The introduction of signed measures in this context is natural, as we recall from the introduction that we study a family of probability measures $\{ \mu_\lambda \}_\lambda$ and their linearization $\mu_{\lambda + \eps} \approx \mu_\lambda + \eps \dot{\mu}$, from which it follows that $\dot{\mu}([0,1]) = 0$, hence $\dot{\mu}$ is in general signed. 

\begin{definition}
    Fix $\lambda \in (0,m)$. We let $T= T_\lambda $ be a linear operator mapping 
    \begin{equs}
    T_\lambda : \CM &\to \CM \\ 
    \eta &\mapsto \sum_{n=4}^\infty n p_n  \int 1_{(\cdot)}\left(\frac{\norm{y}_1  }{\lambda + \norm{y}_1 } \right) \frac{\lambda}{(\lambda + \norm{y}_1 )^2} \dd(\eta \otimes \mu_\lambda^{\otimes (n-1)})(y) \, , 
    \end{equs}
    where $p_n = \P( Z = n)$ and $\norm{y}_1 = y_1 + \dots + y_n $. 
\end{definition}
We note that we can define $T_\lambda$ equivalently by testing $T_\lambda \eta$ against functions $f \in C([0,1])$, resulting in the representation 
\begin{equ}
\int f \dd(T_\lambda \eta) = \sum_{n=4}^\infty np_n \E\left[  \int \frac{\lambda}{(\lambda + x + S_{n-1})^2} f\left( \frac{x + S_{n-1} }{\lambda + x + S_{n-1}} \right) \eta(\dd x) \right]  \, , 
\end{equ}
where $S_{k} = \sum_{i=1}^{k} X_{\lambda, i}$ with $X_{\lambda, i}$ are i.i.d. with law $\mu_\lambda$. We record two immediate properties about $T$.  
\begin{lemma}\label{lem:T_properties}
    \begin{enumerate}
        \item $T_\lambda $ is a positive linear operator, that is $\eta \geq 0  \Rightarrow  T_\lambda  \eta \geq 0$. Furthermore, $\eta \in \CM_{\lambda} \Rightarrow T_\lambda\eta \in \CM_{\lambda}$. 
        \item $T_\lambda$ is bounded in the total variation norm, that is $\norm{T_\lambda \eta}_{\mathrm{TV}} \leq \frac{m}{\lambda} \norm{\eta}_{\mathrm{TV}}$ for any signed measure $\eta \in \CM $. 
    \end{enumerate}
\end{lemma}
\begin{proof}
    The first part follow from straightforward inspection of the definition of $T$; boundedness of the total variation norm follows by the characterization $\norm{\nu}_{\mathrm{TV}} = \sup_{f : |f| \leq 1} \int f \dd \nu $. 
\end{proof}
 Eventually we want to differentiate $A_\lambda$ and $B_\lambda$ with respect to $\lambda$. As $X_\lambda$ is dependent on $\lambda$, the prototypical quantity which we need to differentiate is $\E f(X_\lambda) = \int f(x) \dd\mu_\lambda(x)$. We establish first, that in a suitable way the law $\mu_\lambda $ is differentiable with `derivative' $\Gamma_\lambda$. More precisely: 
\begin{proposition}\label{prop:gamma_fixed_point}
    Let $\lambda \in (0,m)$. 
    \begin{enumerate}
        \item There exists a unique positive measure $\Gamma_\lambda \in \CM_\lambda^+ $ satisfying the fixed point equation
        \begin{equ}\label{eq:gamma_fixed_point}
        \Gamma_\lambda  = \sigma_\lambda + T_\lambda \Gamma_\lambda  \, , 
        \end{equ}
        where $\sigma_\lambda(\dd x) = x(1-x) \mu_\lambda(\dd x)$.  
        \item  The measure $\Gamma_\lambda$ is the derivative of $\mu_\lambda$ in the sense that for all $f \in C^1[0,1]$, the map $\lambda \mapsto \int f \dd \mu_\lambda$ is differentiable and 
        \begin{equ}\label{eq:mu_derivative}
        \frac{\dd}{\dd \lambda} \int f(x) \dd\mu_\lambda(x) = -\frac{1}{\lambda} \int f'(x) \dd \Gamma_\lambda(x) \, . 
        \end{equ}
    \end{enumerate}
\end{proposition}
In order to prove Proposition \ref{prop:gamma_fixed_point} the following bilinear operator will prove useful. 
\begin{definition}
    For fixed $\lambda \in (0,m)$ we define the mapping 
    \begin{equs}
    Q_\lambda :  \CM_\lambda \times \CM_\lambda &\to \R \\ 
    (\eta, \zeta) &\mapsto \int \frac{\lambda}{D_\lambda(x,y)^2} \eta(\dd x) \zeta( \dd y) \, . 
    \end{equs}
\end{definition}
Because of the kernel bounds \eref{eq:D_kernel_bounds} on $D_\lambda$, $Q_\lambda$ is well-defined and it immediately follows that for positive measures $\eta, \zeta \in \CM_\lambda^+$:  
\begin{equ}\label{eq:Q_TV_interplay}
\frac{1}{\lambda} \norm{\eta}_{\mathrm{TV}}\norm{\zeta}_{\mathrm{TV}} \leq  Q_\lambda(\eta, \zeta)  \leq \frac{16}{9\lambda}  \norm{\eta}_{\mathrm{TV}} \norm{\zeta}_{\mathrm{TV}} \, . 
\end{equ}
We will also need a simple technical lemma about the interplay of $Q$ and $T$ and some explicit evaluations of the latter. 
\begin{lemma}\label{lem:interplay_Q_T}
    It holds
    \begin{enumerate}
        \item \label{lem:part1_interplay_Q_T} The linear operator $T_\lambda$ is symmetric with respect to $Q_\lambda$, that is for any signed measures $\eta, \zeta \in \CM_\lambda $ 
        \begin{equ}\label{eq:T_symmetric_to_Q}
        Q_\lambda (T_\lambda \eta, \zeta) = Q_\lambda (\eta, T_\lambda \zeta) \, . 
        \end{equ}
        \item Introducing the measures $\alpha_\lambda, \sigma_\lambda \ll \mu_\lambda $ as 
        \begin{equ}
        \alpha_\lambda(\dd x) = x \mu_\lambda ( \dd x) ,  \, \sigma_\lambda( \dd x) = x(1-x) \mu_\lambda (\dd x)\, ,
        \end{equ}
        they satisfy 
        \begin{equ}\label{eq:alpha_sigma_via_T}
        T_\lambda  \alpha= \sigma \, , 
        \end{equ}
    \end{enumerate}
\end{lemma}
\begin{proof}
    Symmetry of $T$ with respect to $Q$, that is \eref{eq:T_symmetric_to_Q}, follows directly by plugging in the definitions in each other: 
    \begin{equs}
    Q(T \eta, \zeta) &= \sum_{n} np_n \int \lambda \E \left[ \int \frac{\lambda}{(\lambda + z + S_{n-1})^2} D_\lambda\left(\frac{z+S_{n-1}}{\lambda + z + S_{n-1}}, y \right)^{-2} \eta(\dd z) \right] \zeta(\dd y) \\ 
    &=  \sum_{n} np_n \int  \left( \E (\lambda + z + y + S_{n-1} - 1)^{-2} \right) \eta(\dd z) \zeta(\dd y) = Q(\eta, T \zeta) \, .  
    \end{equs}
    For \eref{eq:alpha_sigma_via_T} this follows as straightforwardly, that is, take any $f \in C([0,1])$ and compute: 
    \begin{equs}
    \int f \dd(T \alpha) &= \sum_{n} n p_n \E \int \frac{\lambda y }{(\lambda + y + S_{n-1})^2} f\left(\frac{y+S_{n-1}}{\lambda + y + S_{n-1}}\right) \mu(\dd y) \\ 
    &= \sum_{n} p_n \E \frac{\lambda (nX_n)}{(\lambda + S_n)^2} f\left( \frac{S_n}{\lambda + S_n} \right) = \sum_n p_n \sum_{i=1}^n \E \frac{\lambda X_i }{(\lambda + S_n)^2} f\left( \frac{S_n}{\lambda + S_n} \right)  \\ 
    &= \sum_n p_n \E \frac{\lambda}{(\lambda + S_n)} \frac{S_n}{\lambda + S_n} f\left(\frac{S_n}{\lambda + S_n} \right) = \sum_n p_n \E \left(1- \frac{S_n}{\lambda+ S_n}\right) \frac{S_n}{\lambda + S_n} f\left( \frac{S_n}{\lambda + S_n} \right)  \\ 
    &= \E (1-X)X f(X) = \int f \dd \sigma \, ,  
    \end{equs}
    where in the penultimate equality we used the fixed point of the law of $X$ from \eref{eq:fixed_point_law}. 
\end{proof}
This suffices to proof the first part of Proposition \ref{prop:gamma_fixed_point}. 
\begin{proof}[of 1. of Proposition \ref{prop:gamma_fixed_point}] 
    As usual we split the existence and uniqueness part of the proof.
    \begin{enumerate}
    \item \emph{Existence}: We let $\Gamma^N := \sum_{i=0}^N T^i \sigma$ and aim to show convergence for $N \to \infty$. To this end, first note that as $\sigma \geq 0$, by positivity of $T$ and induction $\Gamma^N \geq 0$. The trick now is to observer that $\alpha - \sigma$ is a positive measure due to  $ \alpha - \sigma = x^2 \mu(\dd x) \geq 0   $ and compute 
    \begin{equs}
     Q(\Gamma^N, x^2 \mu( \dd x)) &= Q(\Gamma^N, \alpha -  \sigma) = \sum_{j=0}^N Q(T^j \sigma, \alpha) - Q(T^j \sigma, \sigma)\\ 
     &= Q(\sigma, \alpha)- Q(T^N \sigma, \sigma) \leq Q(\sigma, \alpha) \, , 
    \end{equs}
    where the penultimate equality comes from $Q(T^j \sigma, \alpha) = Q(T^{j-1}, T\alpha) = Q(T^{j-1}\sigma , \sigma)$ by \eref{eq:T_symmetric_to_Q} and \eref{eq:alpha_sigma_via_T}. Now utilizing the equivalence \eref{eq:Q_TV_interplay} twice, we get: 
    \begin{equ}
    \norm{\Gamma^N}_{\mathrm{TV}} \norm{\alpha - \sigma}_{\mathrm{TV}} \leq \lambda Q(\Gamma^N, \alpha- \sigma) \leq Q(\sigma, \alpha) \leq \frac{16}{9} \norm{\sigma}_{\mathrm{TV}} \norm{\alpha}_{\mathrm{TV}} \, . 
    \end{equ}
    Hence we obtain the uniform in $N$ bound 
    \begin{equ}\label{eq:Gamma_uniform_bound}
    \norm{\Gamma^N}_{\mathrm{TV}} \leq \frac{16}{9} \frac{(\E X_\lambda )  (\E X_\lambda - X_\lambda^2)}{\E X_\lambda^2 } < \infty 
    \end{equ}
    It follows that the sequence of positive measures $(\Gamma^N)_{N \in \N}$ on $[0,1]$ is Cauchy in $(\CM, \norm{\cdot}_{\mathrm{TV}})$ and we denote the limit by $\Gamma := \lim_{N \to \infty} \Gamma^N$. It is immediatley seen that the $\Gamma \in \CM_\lambda^+$, and by continuity of $T$ by Lemma \ref{lem:T_properties} we get 
    \begin{equs}
    T \Gamma = \lim_{N \to \infty} T \Gamma^N = \lim_{N \to \infty} (\Gamma^{N+1} - \sigma) = \Gamma - \sigma \, . 
    \end{equs}
    \item \emph{Uniqueness}: Let $\tilde{\Gamma} \in \CM_\lambda^+ $ another positive measure satisfying the fixed point equation. Then for any $N\geq 1$ we get 
    \begin{equ}
    \tilde{\Gamma} = \underbrace{(\sigma + \dots T^N \sigma )}_{=\Gamma^N } + \underbrace{T^{N+1} \tilde{\Gamma}}_{\geq 0}  \geq \Gamma^N 
    \end{equ}
    Hence taking limits we get $\tilde{\Gamma} \geq \Gamma$. So we apply \eref{eq:Q_TV_interplay}  to $\tilde{\Gamma} - \Gamma , \alpha - \sigma \geq 0 $ and use again symmetry of $T$ with respect to $Q$ and $T\alpha = \sigma$:  
    \begin{equs}
    0 \leq \frac{1}{\lambda} \norm{\tilde{\Gamma} - \Gamma}_{\mathrm{TV}} \norm{\alpha- \sigma}_{\mathrm{TV}} \leq Q(\tilde{\Gamma}- \Gamma , \alpha - \sigma) = Q( T(\tilde{\Gamma} - \Gamma) , \alpha) - Q(\tilde{\Gamma}- \Gamma, \sigma) = 0 \, , 
    \end{equs}
    where the penultimate equality is by the fact that $T(\tilde{\Gamma} - \Gamma) = \tilde{\Gamma} - \Gamma$. Thus we get $\norm{\tilde{\Gamma}- \Gamma }_{\mathrm{TV}} = 0 $. 
\end{enumerate}
\end{proof} 
\begin{remark}
    \begin{enumerate}
        \item Intuitively any $\Gamma$ satisfying \eref{eq:gamma_fixed_point} should be of the form $\Gamma = \sigma + T\sigma + T^2 \sigma + \dots$, so one would be tempted to construct $\Gamma$ via some Neumann series $\sum_k T^k$. Unfortunately, the operator norm of $T$ is only guaranteed to be bounded by  $ \frac{m}{\lambda} \geq 1$, which is the reason for the above proof. 
        \item One can attach a probabilistic meaning to $\Gamma_\lambda$ in the following way. We formally differentiate \eref{eq:fixed_point_law} to get 
        \begin{equs}
        \partial_\lambda X_\lambda &= - \frac{S_\lambda}{(\lambda + S_\lambda)^2} + \frac{\lambda}{(\lambda + S_\lambda)^2} \sum_i X_{\lambda, i} \\ 
        &= - \frac{1}{\lambda} X_\lambda (1-X_\lambda) + X_\lambda (1-X_\lambda) \frac{\sum_{i=1}^Z X_{\lambda, i}}{S_\lambda}  \, . 
        \end{equs}
        Rewriting this further, while emphasizing the dependence on the vertex $v$ yields 
        \begin{equ}\label{eq:stopping_time_recursion_prelim}
        -\lambda \partial_\lambda \log X_\lambda(v)  = (1-X_\lambda(v))\left[ 1 + \frac{\sum_{i=1}^Z (-\lambda \partial_\lambda \log X_{\lambda}(v_i)) X_{\lambda}(v_i) }{S_\lambda}  \right]  \, , 
        \end{equ}
        where $v_i$ is the $i$-th child of $v$.
        It turns out that \eref{eq:stopping_time_recursion_prelim} can be given a meaning in terms of a stopping time $\tau(v)$ of an auxiliary random walk on the given Galton Watson tree $\CT$. This auxiliary random walk is defined via transition probabilities at vertex $v$ as 
        \begin{enumerate}
            \item with probability $X_\lambda(v)$ it becomes extinct. 
            \item with probability $1-X_\lambda(v)$ it moves to one of the children of $v_i$, each with probability $\frac{X_\lambda(v_i)}{S_\lambda}$. 
        \end{enumerate}
        With respect to this construction, let $\tau(v)$ be the stopping time until extinction of a random walk started at the vertex $v$. We can immediately see by a first step analysis that 
        \begin{equ}
        \E [ \tau(v) | \CT ] = (1-X_\lambda)\left(1 + \sum_{i=1}^Z \frac{X_\lambda(v_i)}{S_\lambda(v)} \E [ \tau(v_i) | \CT] \right) \, . 
        \end{equ}
        Hence this suggests, assuming uniqueness of solutions, that  
        \begin{equ}\label{eq:stopping_time_equality}
        \E [ \tau(v) | \CT ] = -\lambda \partial_\lambda \log X_\lambda(v) \, . 
        \end{equ}
        To connect this now to a representation of $\Gamma_\lambda$, we furthermore make a formal calculation 
        \begin{equ}
        - \frac{1}{\lambda} \int f'(x) \Gamma_\lambda(\dd x)= \partial_\lambda \E f(X_\lambda) = \E [  f'(X_\lambda) \partial_\lambda X_\lambda ]  = \int f'(x) \E \left[ \partial_\lambda X_\lambda \delta_{X_\lambda}(\dd x) \right]  \, . 
        \end{equ}
        from which we are able to identify 
        \begin{equ}
        \Gamma_\lambda(\dd x) = - \lambda \E [ \partial_\lambda X_\lambda \delta_{X_\lambda}(\dd x) ] \, . 
        \end{equ}
        Substituting \eref{eq:stopping_time_equality} we arrive at: 
        \begin{equ}\label{eq:informal_representation_gamma}
        \Gamma_\lambda(\dd x) = \E [ X_\lambda \tau(v)  \delta_{X_\lambda}( \dd x) ] \, . 
        \end{equ}
        \item Starting from \eref{eq:informal_representation_gamma}, we can also make sense of the Dyson series in a probabilistic manner. To this end note that by the layer-cake representation 
        \begin{equ}
        \Gamma_\lambda(\dd x) = \E X_\lambda \tau(\rho) \delta_{X_\lambda}(\dd x) = \sum_{n=1}^\infty \E X_\lambda 1_{\tau \geq n} \delta_{X_\lambda}(\dd x) = \sum_{n=1}^\infty \E \left[ X_\lambda \E[ 1_{\tau \geq n} | \CT ] \delta_{X_\lambda}(\dd x) \right]  \, . 
        \end{equ}
        Let us compute now explicitly 
        \begin{enumerate}
            \item $n=1$: It is immediate that by definition of the auxiliary process $\P( \tau(\rho) \geq 1| \CT) = (1-X_\lambda (\rho))$. Hence 
            \begin{equ}
            \E \left[ X_\lambda \E[ 1_{\tau \geq 1} | \CT ] \delta_{X_\lambda}(\dd x) \right]  = \E[ X_\lambda (1-X_\lambda) \delta_{X_\lambda}(\dd x)] = x(1-x)\mu_\lambda(\dd x) = \sigma_\lambda(\dd x) \, . 
            \end{equ}
            \item $n=2$: Again, by inspection of the auxiliary random walk, we see the probability that the walk doesn't get extinct for at least two steps is 
            \begin{equ}
            \P(\tau \geq 2| \CT) = (1-X_\lambda(\rho))\left[ \sum_i \frac{X_\lambda(\rho_i)}{S_\lambda(\rho)}(1-X_\lambda(\rho_i))  \right]  \, . 
            \end{equ}
            And so 
            \begin{equs}
            \E \left[ X_\lambda \E[ 1_{\tau \geq 2} | \CT ] \delta_{X_\lambda}(\dd x) \right]  &= \E \left[  \frac{\lambda}{(\lambda + S_\lambda(\rho))^2} \left( \sum_{i} X_\lambda(\rho_i)(1- X_\lambda(\rho_i))  \right) \delta_{X_\lambda}(\dd x)  \right]  \, . 
            \end{equs}
            Testing against some function $f$, we see 
            \begin{equs}
            \int f(x) \E \left[ X_\lambda \E[ 1_{\tau \geq 2} | \CT ] \delta_{X_\lambda}(\dd x) \right] &= \sum_n np_n \int \frac{\lambda}{(\lambda + \norm{x}_1)^2} f\left( \frac{\norm{x}_1}{\lambda + \norm{x}_1} \right) \dd(\mu_\lambda^{n-1}\otimes \sigma_\lambda)(x) \\ 
             &= \int f \dd(T_\lambda \sigma_\lambda) \, . 
            \end{equs}
        \end{enumerate}
        And generally for $n \in \N$ 
        \begin{equ}
         \E \left[ X_\lambda \E[ 1_{\tau \geq n} | \CT ] \delta_{X_\lambda}(\dd x) \right] = T_\lambda^{n-1} \sigma_\lambda (\dd x) \, . 
        \end{equ}
        Thus we have rederived the Dyson series $\Gamma = \sigma + T\sigma + \dots$ in a purely probabilistic manner. It is also worthwhile to note that the analytic convergence question of $\Gamma^N$ can be answered (with a bit more work to make the previous steps fully rigorous) by the probabilistic statement of (the weighted) integrability of the existence time $\tau$ of the auxiliary process : 
        \begin{equ}
        \norm{\Gamma_\lambda}_{\mathrm{TV}} = \E \left[ X_\lambda(\rho) \tau(\rho)  \right] < \infty \, . 
        \end{equ}
    \end{enumerate}
\end{remark}
To prove the second part, we will need a technical lemma showing a certain Lipschitz continuity in law: 
\begin{lemma}\label{lem:X_lipschitz}
    For any $u,v \geq 0$ such that $[u,v] \subset (0,m)$ there exists a constant $C_{u,v} > 0 $ such that 
    \begin{equ}
    \E \left[ X_a - X_b  \right] \leq C_{u,v} (b-a) 
    \end{equ}
    uniformly in $a,b$ with $u \leq a < b \leq v$. 
\end{lemma}
\begin{proof}
    See Lemma \ref{lem:X_lipschitz_appendix}. 
\end{proof}
\begin{proof}[of 2. of Proposition \ref{prop:gamma_fixed_point}]
    Let $f \in C^1([0,1])$.  For $a,b \in (0,m)$ with $a<b$ one gets using $X_a \geq X_b $ (pathwise for every fixed tree by an application of Rayleigh's monotonicity principle \cite[Section 2.4]{lyons2005probability}), that 
    \begin{equs}
    f(X_b) - f(X_a) = - \int_{[0,1]} 1_{[X_b, X_a]}(y) f'(y) \dd y \, , 
    \end{equs}
    and hence by taking the expectation we get the finite difference equation: 
    \begin{equ}
    \frac{\int f \dd \mu_b - \int f \dd \mu_a}{b-a} = -\frac{1}{a} \left( \int f'(y) \underbrace{\frac{a}{b-a}\E \left[1_{[X_b, X_a](y)} \right] \dd y}_{=: R_{a,b}(\dd y)}   \right) = -\frac{1}{a} \int f'(y) R_{a,b}(\dd y) \, . 
    \end{equ}
    The rest is now rigorously establishing that we can take the limit as $b \to a$ in the following steps. 
    \begin{enumerate}
        \item The map $\lambda \mapsto \E f(X_\lambda), (0,m) \to \R $ is locally Lipschitz: Take some $u,v$ with $[u,v] \subset (0,m)$, then for any $a,b$ with $u \leq a < b \leq v$ we have 
        \begin{equs}
        | \E f(X_b) - \E f(X_a) | &\leq \mathrm{Lip}(f) \E |X_b - X_a| \leq \mathrm{Lip}(f) C_{u,v} (b-a)  \, , 
        \end{equs}
        by Lemma \ref{lem:X_lipschitz}. Hence, by Rademacher's theorem, for Lebesgue almost every $\lambda \in (0,m)$ the map $\lambda \mapsto \E f(X_\lambda)$ is differentiable. By the usual arguments, for Lebesgue almost every $a \in (0,m)$, for every rational coefficient polynomial $f$ the following exists 
        \begin{equ}
        \frac{\dd }{\dd \lambda } \E f(X_\lambda ) \big|_{\lambda = a} \in \R \, . 
        \end{equ}
        For further readability, call $\CP_\Q \subset C([0,1])$ the countable set of rational coefficient polynomials and $\CN \subset (0,m) $ the exceptional set of Lebesgue measure zero. 
        \item Take any $a \in (0,m) \setminus \CN$  and choose a sequence $b_n \downarrow a$. Looking at the positive measure $R_{a,b}$, we get immediatley by Lemma \ref{lem:X_lipschitz}
        \begin{equ}
        \norm{R_{a,b_n}}_{\mathrm{TV}} = R_{a,b_n}([0,1]) = a \frac{\E X_a - X_{b_n} }{b_n-a} \leq a C_{a , b_0} 
        \end{equ}
        Hence the family $(R_{a,b_n})_n \subset \CM^+$ is tight, and so we can extract a subsequential weak limit which we call 
        \begin{equ}
        \tilde{\Gamma}_a := \weaklim_{k \to \infty } R_{a,b_{n_k}} \, . 
        \end{equ}
        \item We can now identify 
        \begin{equs}
        -\frac{1}{a} \int f'(y) \tilde{\Gamma}_a(\dd y) = -\frac{1}{a} \lim_{k \to \infty} \int f'(y) R_{a,b_{n_k}}(\dd y) = \lim_{k \to \infty} \frac{\E X_{b_{n_k}} - \E X_a}{b-a} = \frac{\dd }{\dd \lambda } \E f(X_\lambda ) \big|_{\lambda = a}  \, .  
        \end{equs}
        \item Moving on to show the fixed point equation, we again note that $X_\lambda = \frac{\sum_{i=1}^Z X_{\lambda, i}}{\lambda + \sum_{i=1}^Z X_{\lambda, i}}$ in distribution, that is, for any $\lambda \in (0,m)$ and  $f \in \CP_\Q $ we have 
        \begin{equ}
        \int f(x) \dd\mu_\lambda(x) = \E f(X_\lambda ) = \E f\left( \frac{S_\lambda}{\lambda + S_\lambda }\right) =  \sum_{n=4}^\infty p_n  \int \phi_{\lambda, n}  \dd \mu^{\otimes n}_\lambda \, . 
        \end{equ}
        where $p_n = \P(Z =n)$ and $\phi_{\lambda,n} : y \mapsto f\left(\frac{\norm{y}_1  }{\lambda + \norm{y}_1 }\right)$. Here we can write again $\norm{y}_1 = y_1 + \dots + y_n$ since $y \in [0,1]^n$. We can differentiate this fixed point equation at $\lambda \in (0,m) \setminus \CN $ to get 
        \begin{equs}
        - \frac{1}{\lambda} \int f'(x) \dd \Gamma_\lambda(x) &= \sum_n p_n \left[  \int \partial_\lambda \phi_{\lambda, n} \dd \mu_\lambda^{\otimes n}   - \frac{n}{\lambda} \int  \partial_1 \phi_{\lambda, n} \dd \left( \tilde{\Gamma}_\lambda \otimes \mu_\lambda^{\otimes (n-1)} \right)  \right]  \, . 
        \end{equs}
        Here, to arrive rigorously at the right-hand-side we employed Lemma \ref{lem:rigorous_derivative_fixed_point}. Now, since 
        \begin{equs}
        \partial_\lambda \phi_{\lambda, n}(y ) &= f'\left(  \frac{\norm{y}_1  }{\lambda +\norm{y}_1  }\right) \left( -\frac{\norm{y}_1 }{(\lambda + \norm{y}_1 )^2} \right) \\ 
        \partial_1 \phi_{\lambda, n}(y) &= f'\left(  \frac{\norm{y}_1  }{\lambda + \norm{y}_1  }\right) \frac{\lambda}{(\lambda + \norm{y}_1 )^2} \, , 
        \end{equs}
        the first term can be seen to be 
        \begin{equ}
        \int \partial_\lambda \phi_{\lambda, n} \dd \mu_\lambda^{\otimes n} = - \frac{1}{\lambda}  \int f'\left(\frac{\norm{y}_1  }{\lambda + \norm{y}_1 }\right)\frac{\norm{y}_1  }{\lambda + \norm{y}_1}\left(1-\frac{\norm{y}_1 }{\lambda + \norm{y}_1 }\right) \dd \mu_\lambda^{\otimes n}(y)   = -  \frac{1}{\lambda} \E \left[ f'(X) X(1-X) | Z= n \right] \, , 
        \end{equ}
        and hence summing over $n$ yields directly $\int f'(x) \dd \sigma_\lambda(x) $. 
        The second term is already in the form of the operator $T$ and so we are done. Hence we arrive at  
        \begin{equ}
        \int f' \dd \tilde{\Gamma}_\lambda = \int f' \dd \sigma_\lambda + \int f' \dd(T \tilde{\Gamma}_\lambda) \, , 
        \end{equ}
        and as $f \in \CP_\Q $ and the polynomials are dense in $C([0,1])$, we get $\tilde{\Gamma}_\lambda = \sigma_\lambda + T \tilde{\Gamma}_\lambda$. By uniqueness of the fixed point equation \eref{eq:gamma_fixed_point}, we get $\tilde{\Gamma}_a = \Gamma_a $ for any $a \in (0,m) \setminus \CN$. 
        \item Finally we want to upgrade the convergence result. So fix $f' \in \CP_\Q$. From the previous steps we know that $\lambda \mapsto \E f(X_\lambda)$ is locally Lipschitz continuous and so on every compact interval $[a,b] \subset (0,m)$ it is absolutely continous. As such, by Lebesgue's fundamental theorem of calculus, for any $a,b \in (0,m)$: 
        \begin{equ}
        \E f(X_b) - \E f(X_a) = \int_{[a,b]} \frac{\dd}{\dd \lambda} \E f(X_\lambda) \dd \lambda = \int_{[a,b]} \left( -\frac{1}{\lambda} \int f'(y) \Gamma_\lambda( \dd y)  \right) \dd \lambda \, , 
        \end{equ}
        where the last equality comes from the almost everywhere equivalence. Using the fact that $\lambda \mapsto \Gamma_\lambda$ is weakly continuous by Lemma \ref{lem:weak_continuity}, we know that $\lambda \mapsto -\frac{1}{\lambda} \int f'(y) \Gamma_\lambda(\dd y) $ is continuous, and hence the right-hand-side is equivalent to the Riemann integral. Hence, by the fundamental theorem of calculus 
        \begin{equ}
        \forall \lambda \in (0,m): \quad \frac{\dd}{\dd \lambda} \E f(X_\lambda)  = -\frac{1}{\lambda} \int f'(y) \Gamma_\lambda(\dd y) \, . 
        \end{equ}
        Now it is straightforward to extend the argument to $f \in  C^1([0,1]) \subset  \overline{\CP_\Q}^{\norm{\cdot}_\infty} $. 
    \end{enumerate}
\end{proof}
Having established the derivative formula \eref{eq:mu_derivative}, we record in passing that it is not difficult to show (using weak continuity of $\lambda \mapsto \mu_\lambda$ from Lemma \ref{lem:weak_continuity}) that for $f: (0,m) \times [0,1]^n \to \R$ continuously differentiable, we have the formula 
\begin{equ}\label{eq:derivative_lambda_iterated_integral}
\frac{\dd}{\dd \lambda} \int f_\lambda \dd\mu_\lambda^{\otimes n} = \int \partial_\lambda f_\lambda \dd \mu_\lambda^{\otimes n} - \frac{1}{\lambda} \sum_{i=1}^n \int \partial_i f_\lambda \dd (\mu_\lambda^{\otimes(i-1)} \otimes \Gamma_\lambda \otimes \mu_\lambda^{\otimes(n-i)}) \, . 
\end{equ}
We can now explicitly compute the derivative of \eref{eq:rousselin_formula}. 
\begin{proposition}\label{prop:ell_prime_formula}
    The speed $\ell$ is differentiable on $(0,m)$ and it holds that 
    \begin{equ}\label{eq:ell_prime_formula}
    -\lambda^2 B_\lambda \ell'(\lambda) = 2 Q_\lambda(\Omega_\lambda, w_\lambda) + (1+\ell(\lambda) ) Q_\lambda(\sigma_\lambda, \sigma_\lambda)\, . 
    \end{equ}
    where  $\Omega_\lambda = \sigma_\lambda + (\lambda - 1) \Gamma_\lambda$ and $ w_\lambda(\dd x) = (x - \ell(\lambda)(1-x))\mu_\lambda(\dd x) $.    
\end{proposition}
\begin{proof}
    We start with the rewriting of Rousselins formula \eref{eq:rousselin_formula}
    \begin{equ}
    0 = A_\lambda - B_\lambda \ell(\lambda) = \E \frac{XY - \ell(\lambda)(X+Y- XY)}{D_\lambda(X,Y)} =: \E \frac{H_\lambda(X,Y)}{D_\lambda(X,Y)} = \int \frac{H_\lambda(x,y)}{D_\lambda(x,y)} \dd\mu_\lambda(x) \dd \mu_\lambda(y) \, , 
    \end{equ}  and differentiate with respect to $\lambda$. We obtain with \eref{eq:derivative_lambda_iterated_integral} \footnote{even though $H_\lambda / D_\lambda$ is not $C^1$ on $(0,m) \times [0,1]^2$, we can easily remedy this fact by noting that $\mu_\lambda$ has restricted support on which $1/D_\lambda$ \emph{is} $C^1$, so the usual argument with smooth bump functions does the job. }
    \begin{equ}
    0 = \int \partial_\lambda \left( \frac{H_\lambda(x,y)}{D_\lambda(x,y)} \right)  \dd\mu_\lambda(x) \dd \mu_\lambda(y) - \frac{1}{\lambda} \partial_x \left( \frac{H_\lambda(x,y)}{D_\lambda(x,y)} \right)  \dd\Gamma_\lambda(x) \dd \mu_\lambda(y) - \frac{1}{\lambda} \partial_y \left( \frac{H_\lambda(x,y)}{D_\lambda(x,y)} \right)  \dd\Gamma_\lambda(y) \dd \mu_\lambda(x) \, . 
    \end{equ}
    Performing the algebraic manipulations and correctly computing the derivatives leads to 
    \begin{equ}
    -\lambda^2 B_\lambda \ell'(\lambda) = \lambda^2 \E \left[\frac{H_\lambda(X,Y)}{D^2_\lambda(X,Y)} \right] + 2 Q(\Gamma, x^2 \mu_\lambda(\dd x) ) + 2(\lambda  - 1) Q(\Gamma, w) \, . 
    \end{equ}
    Here $w(\dd x)= (x - \ell(\lambda)(1-x)) \mu_\lambda(\dd x) $. We continue with the rewriting by noting that $\lambda = D_\lambda + (1-x)(1-y)$ and hence 
    \begin{equ}
     \lambda^2 \E \left[\frac{H_\lambda(X,Y)}{D^2_\lambda(X,Y)} \right] = \lambda \E \left[ \frac{(D_\lambda(X,Y) + (1-X)(1-Y))H_\lambda(X,Y)}{D_\lambda^2(X,Y)} \right] = \lambda \E \frac{(1-X)(1-Y)H_\lambda(X,Y)}{D_\lambda^2(X,Y)} \, , 
    \end{equ}
    where the last equality follows as we have $\E \frac{H}{D} = 0 $. By noticing that 
    \begin{equ}
    (1-x)(1-y)H_\lambda(x,y) = (1+\ell)(x(1-x))(y(1-y)) - \ell(x(1-x)(1-y) + y(1-y)(1-x)) \, , 
    \end{equ}
    we arrive at 
    \begin{equ}
    -\lambda^2 B_\lambda \ell'(\lambda) = (1+\ell) Q(\sigma, \sigma) - 2 \ell(\lambda) Q(\sigma, (1-x) \mu_\lambda(\dd x)) + 2Q(\Gamma, x^2 \mu_\lambda(\dd x)) + 2(\lambda - 1) Q(\Gamma, w) \, . 
    \end{equ}
    Now we are almost there, we use that $\Gamma = \sigma + T\Gamma$ and that $\sigma = T\alpha$ to get at 
    \begin{equs}
    Q(\Gamma, x^2 \mu_\lambda(\dd x)) &= Q(\Gamma, x \mu_\lambda(\dd x)) - Q(\Gamma, \sigma) = Q(\Gamma, x \mu_\lambda(\dd x)) - Q(\Gamma, T (x \mu_\lambda(\dd x))) \\ 
    &= Q(\Gamma - T \Gamma , x \mu_\lambda(\dd x)) = Q(\sigma, x \mu_\lambda(\dd x)) \, . 
    \end{equs}
    And hence 
    \begin{equ}
    -\lambda^2 B_\lambda \ell'(\lambda) = (1+\ell)Q(\sigma,\sigma) + 2 Q(\sigma + (\lambda-1)\Gamma, w) \, . 
    \end{equ}
\end{proof}
\section{The proof}
Having established \eref{eq:ell_prime_formula}, we can now quite straightforwardly prove Theorem \ref{thm:lpp_four_children}. We start with a technical bound on $\ell(\lambda)$ which is independent of the assumption on the offspring law of $Z \geq 4$. 
\begin{lemma}
    For $\lambda \in (0,m)$ it holds that 
    \begin{equ}\label{eq:ell_upper_bound} 
    \ell(\lambda) \leq \frac{\E X_\lambda}{2- \E X_\lambda} \, . 
    \end{equ}
\end{lemma}
\begin{proof}
    The observation we use is that 
    \begin{equ}\label{eq:A_lambda_ineq}
    A_\lambda = \E \frac{XY}{D_\lambda(X,Y)} \leq \E \left[ Y \right] \E \left[ \frac{X}{D_\lambda(X,Y)} \right]   \, . 
    \end{equ}
    Indeed, the function $r: y \mapsto \E \frac{X}{D_\lambda(X,y)} = \E \frac{X}{\lambda - (1-X)(1-y)}$ is monotonically decreasing, hence we have for any random variables $Y,Z$  that 
    \begin{equ}
    (Y - Z )(r(Y) - r(Z)) \leq 0 \, . 
    \end{equ}
    Hence taking in particular $Y, Z$ i.i.d. and distributed according to $\mu_\lambda$, we get by taking expectation values the desired inequality \eref{eq:A_lambda_ineq}. Finally, by definition we have 
    \begin{equ}
    B_\lambda = \E \frac{X+Y - XY}{D_\lambda(X,Y)} = 2 \E\left[ \frac{X}{D_\lambda(X,Y)} \right] - A_\lambda \, , 
    \end{equ}
    and so 
    \begin{equ}
    \ell(\lambda) = \frac{A_\lambda}{B_\lambda}= \frac{A_\lambda}{2\E\left[ \frac{X}{D_\lambda(X,Y)} \right] - A_\lambda } \leq \frac{ \E \left[ Y \right] \E \left[ \frac{X}{D_\lambda(X,Y)} \right] }{ 2 \E\left[ \frac{X}{D_\lambda(X,Y)} \right]   -  \E \left[ Y \right] \E \left[ \frac{X}{D_\lambda(X,Y)} \right] } = \frac{\E Y}{2- \E Y} \, . 
    \end{equ}
\end{proof}
Having established now this bound on the speed $\ell$ we come now to the proof of Theorem \ref{thm:lpp_four_children}. The point is that we can lower-bound very crudely each of the appearing terms inside of \eref{eq:ell_prime_formula}:
\begin{proof}[of Theorem \ref{thm:lpp_four_children}]
    We start by inspecting the formula \eref{eq:ell_prime_formula} and note that as $Q_\lambda(\sigma_\lambda, \sigma_\lambda) \geq 0$ we are done the moment we can show that the sign of $Q_\lambda(\Omega_\lambda, w_\lambda)$ is (strictly) positive as well. Explicitly we can write  
    \begin{equs}
    Q_\lambda(\Omega_\lambda, w_\lambda) &= \int \int \frac{y - \ell(\lambda)(1-y)}{D^2_\lambda(x,y)} \mu_\lambda(\dd y) \Omega_\lambda(\dd x)   
    \end{equs}
    and the inner integrand we can estimate directly for every fixed  $x$ 
    \begin{equs}
    \int \frac{y- \ell(\lambda)(1-y)}{D_\lambda^2(x,y)} \mu_\lambda(\dd y) &\geq \int \frac{y }{\lambda^2} \mu_\lambda(\dd y) - \ell(\lambda) \int \frac{1-y}{(\frac{3}{4}\lambda)^2} \mu_\lambda(\dd y)  \\ 
    &\geq  \frac{\E X_\lambda}{\lambda^2} - \frac{\E X_\lambda}{2-\E X_\lambda} \frac{16}{9\lambda^2}  (1-\E X_\lambda) \\ 
    &= \frac{\E X_\lambda}{\lambda^2} \left( 1- \frac{16}{9} \frac{1-\E X_\lambda}{2-\E X_\lambda} \right) \geq \frac{\E X_\lambda}{9 \lambda^2}\, > 0 .  
    \end{equs}
Here in the second inequality we used \eref{eq:ell_upper_bound} and in the last inequality the elementary inequality $\frac{1-x}{2-x} \geq \frac{1}{2}$. Since $\Omega_\lambda$ is a positive measure by Lemma \ref{lem:omega_positive_measure}, we obtain that $Q_\lambda(\Omega_\lambda, w_\lambda) > 0 $ and hence that $-\lambda^2 B_\lambda \ell'(\lambda) > 0 $ and so the theorem follows, as $\lambda^2 B_\lambda > 0 $. 
\end{proof}

\appendix 
\section{Some technical lemmas}
We collect some lemmas that were used, which are technical in nature and use often just a slight (but tedious) variant of arguments already presented in the main text.
\begin{lemma}\label{lem:X_lipschitz_appendix}
    For any $u,v \geq 0$ such that $[u,v] \subset (0,m)$ there exists a constant $C_{u,v} > 0 $ such that 
    \begin{equ}
    \E \left[ X_a - X_b  \right] \leq C_{u,v} (b-a) 
    \end{equ}
    uniformly in $a,b$ with $u \leq a < b \leq v$. 
\end{lemma}
\begin{proof}
    We fix such $u,v$ as in the assumption and choose some $a,b$ with $u \leq a < b \leq v$. We will now define depending on $a,b$ analogous objects to $Q,T$ and  $\alpha, \sigma$ as follows. First we set 
    \begin{equ}
    D_a(x,x') = a - (1-x)(1-x'), \quad D_b(y,y') = b - (1-y)(1-y') \, . 
    \end{equ} 
    The symmetric bilinear form $Q_{a,b}$ we set as 
    \begin{equs}
    Q_{a,b} : \CM_{a,2} \times \CM_{b,2} &\to \R \\ 
        (\eta, \zeta) &\mapsto \int \frac{\sqrt{ab}}{D_a(x,x') D_b(y,y')} \eta(\dd x, \dd x') \zeta(\dd y, \dd y') \, . 
    \end{equs}
    Here, $\CM_{\lambda,2} = \{ \nu | \supp \nu \subset [( 1- \lambda / 4)_+, 1]^2 \} $. We get similar to \eref{eq:Q_TV_interplay}, that for positive measures $\eta \in \CM_{a,2}^+, \zeta \in \CM_{b,2}^+  $:  
    \begin{equ}\label{eq:Q_ab_TV_interplay}
    \frac{1}{v} \norm{\eta}_{\mathrm{TV}} \norm{\zeta}_{\mathrm{TV}} \leq Q_{a,b}(\eta, \zeta) \leq \frac{16}{9u} \norm{\eta}_{\mathrm{TV}} \norm{\zeta}_{\mathrm{TV}} \, . 
    \end{equ}
    as well as an operator $T_{a,b} : \CM \to \CM$ defined via actions on test functions $f \in C([0,1]^2)$ as 
    \begin{equ}
    \int f \dd (T_{a,b} \eta) = \sum_{n} np_n \int \E \left[ \frac{\sqrt{ab}}{(a+x+ S_{a,n-1})(b+ y + S_{b,n-1})} f\left(\frac{x+S_{a,n-1}}{\lambda + x + S_{a,n-1}},\frac{y+S_{b,n-1}}{\lambda + y + S_{b,n-1}} \right) \right] \eta(\dd x , \dd y) \, . 
    \end{equ}
    This map $T_{a,b}$ has analogous properties as in Lemmas \ref{lem:T_properties} and \ref{lem:interplay_Q_T}, in particular it is a positive linear operator symmetric with respect to $Q_{a,b}$, that is for $\eta \in \CM_{a,2}, \zeta \in \CM_{b,2}$ 
    \begin{equ}
    Q_{a,b}(T_{a,b}\eta, \zeta) = Q_{a,b}(\eta, T_{a,b} \zeta) \, . 
    \end{equ}
    We introduce furthermore 
    \begin{equ}
    \mu_{a,b} = \mathrm{Law}_\P(X_a, X_b) \, , 
    \end{equ}
    as well as 
    \begin{equs}
    \alpha_{a,b}(\dd x, \dd y) &= \left( \frac{x}{\sqrt{a}} + \frac{y}{\sqrt{b}} \right) \mu_{a,b}(\dd x, \dd y) \\ 
    \sigma_{a,b}(\dd x, \dd y) &= \left( \frac{x(1-y)}{\sqrt{b}} + \frac{y(1-x)}{\sqrt{a}} \right) \mu_{a,b}(\dd x , \dd y) \, . 
    \end{equs}
    We can check by similar arguments as before that 
    \begin{equs}
    T_{a,b} \alpha_{a,b} &= \sigma_{a,b} \, . 
    \end{equs}
    Finally we define 
    \begin{equ}
    \Gamma_{a,b}( \dd x , \dd y) = \frac{x-y}{\sqrt{b} - \sqrt{a}} \mu_{a,b}(\dd x, \dd y) \, . 
    \end{equ}
    We can now by straightforward algebraic manipulations check that in fact 
    \begin{equ}
    \Gamma_{a,b} = \sigma_{a,b} + T_{a,b} \Gamma_{a,b} \, . 
    \end{equ}
    Now we have done enough preprocessing to obtain a bound on $\norm{\Gamma_{a,b}}_{\mathrm{TV}}$. In fact, the procedure will be mimicking again the main text. That is, we note that $\alpha_{a,b} - \sigma_{a,b}$ is positive: 
    \begin{equ}
    \frac{\dd (\alpha_{a,b} - \sigma_{a,b}) }{\dd \mu_{a,b}}(x,y)  = \left((x-y)\left( \frac{1}{\sqrt{a}} - \frac{1}{\sqrt{b}} \right) + xy \left( \frac{1}{\sqrt{a}} + \frac{1}{\sqrt{b}} \right)   \right)   \geq 0 \, . 
    \end{equ}
    Then we compute using symmetry of $T_{a,b}$ with respect to $Q_{a,b}$ and $T_{a,b}\alpha_{a,b} = \sigma_{a,b}$: 
    \begin{equs}
    Q_{a,b}(\Gamma_{a,b}, \alpha_{a,b} - \sigma_{a,b}) = Q_{a,b}(\sigma_{a,b}, \alpha_{a,b}) \, . 
    \end{equs}
    Hence, with \eref{eq:Q_ab_TV_interplay} we obtain 
    \begin{equ}
    \frac{1}{v} \norm{\Gamma_{a,b}}_{\mathrm{TV}} \norm{\alpha_{a,b} - \sigma_{a,b}}_{\mathrm{TV}} \leq Q_{a,b}(\Gamma_{a,b}, \alpha_{a,b} - \sigma_{a,b}) = Q_{a,b}(\sigma_{a,b}, \alpha_{a,b}) \leq \frac{16}{9u} \norm{\alpha_{a,b}}_{\mathrm{TV}} \norm{\sigma_{a,b}}_{\mathrm{TV}} \, . 
    \end{equ}
    Now all that is left is to explicitly compute bounds on the total mass of $\alpha_{a,b}, \sigma_{a,b} $ and their difference: 
    \begin{equs}
    \norm{\alpha_{a,b}}_{\mathrm{TV}} &= \frac{1}{\sqrt{a}} \E X_a + \frac{1}{\sqrt{b}} \E X_b \leq \frac{2}{\sqrt{u}}  \\ 
    \norm{\sigma_{a,b}}_{\mathrm{TV}} &= \frac{1}{\sqrt{b}} \E X_a(1-X_b) + \frac{1}{\sqrt{b}} \E X_b(1-X_a) \leq \frac{2}{\sqrt{u}} 
    \end{equs}
    For the difference, we want to compute a lower bound 
    \begin{equ}
     \norm{\alpha_{a,b} - \sigma_{a,b} }_{\mathrm{TV}} \geq \left( \frac{1}{\sqrt{a}} + \frac{1}{\sqrt{b}} \right) \E X_a X_b \geq \frac{2}{\sqrt{v}} \E X_v^2 > 0 \, . 
    \end{equ}
    Consequently we obtain 
    \begin{equ}
    \norm{\Gamma_{a,b}}_{\mathrm{TV}} \leq \frac{32 v^{3/2}}{9 u^2 \E X_v^2} \, . 
    \end{equ}
    Hence, we get the desired Lipschitz in law estimate: 
    \begin{equs}
    \E \left[ X_a - X_b \right]  = (\sqrt{b} - \sqrt{a}) \norm{\Gamma_{a,b}}_{\mathrm{TV}} \leq \underbrace{\left(  \frac{16 v^{3/2}}{9 u^{5/2} \E X_v^2}  \right)}_{=: C_{u,v}} (b-a) \, . 
    \end{equs}
\end{proof}

\begin{lemma}\label{lem:rigorous_derivative_fixed_point}
    In the context of the proof for Proposition \ref{prop:ell_prime_formula}, it holds that for $\lambda \in (0,m) \setminus \CN$ once we set $\phi_{\lambda,n}(y) = f\left( \frac{y_1 + \dots + y_n}{\lambda + y_1 + \dots + y_n} \right)$ we have 
    \begin{equ}
    \frac{\dd }{\dd \lambda} \sum_n p_n \int \phi_{\lambda, n} \dd \mu_\lambda^{\otimes n} = \sum_n p_n \left[ \int \partial_\lambda \phi_{\lambda, n} \dd \mu_\lambda^{\otimes n} - \frac{n}{\lambda} \int   \partial_1\phi_{\lambda, n} \dd\left( \tilde{\Gamma}_\lambda \otimes \mu_\lambda^{\otimes(n-1)} \right)  \right] \, .
    \end{equ}
\end{lemma}
\begin{proof}
    We evaluate the derivative by taking the limit of the finite difference quotient along a specific sequence $\lambda_k \downarrow  \lambda$ for which we know that $R_{\lambda, \lambda_k} \xrightharpoonup{} \tilde{\Gamma}_\lambda$. Without loss of generality this sequence will be contained in some compact interval $(\lambda_k)_k \subset [u,v] \subset (0,m)$.  That is we form 
    \begin{equ}
    \frac{\E f(X_{\lambda_k}) - \E f(X_\lambda)}{\lambda_k - \lambda} = \sum_n p_n \frac{\int \phi_{\lambda_k, n} \dd \mu_{\lambda_k}^{\otimes n} - \int \phi_{\lambda, n} \dd \mu_\lambda^{\otimes n}}{\lambda_k - \lambda } =: \sum_n p_n I_n(\lambda, \lambda_k) 
    \end{equ}
    We split now $I_n = I_n^{(1)} + I_n^{(2)}$ with 
    \begin{equs}
     I_n^{(1)}(\lambda, \lambda_k ) = \frac{\int \left( \phi_{\lambda_k,n} - \phi_{\lambda, n} \right)  \dd \mu_{\lambda_k}^{\otimes n} }{\lambda_k - \lambda}, \quad  
     I_n^{(2)}(\lambda, \lambda_k) = \frac{\int \phi_{\lambda, n} \dd (\mu_{\lambda_k}^{\otimes n} - \mu_{\lambda}^{\otimes n})}{\lambda_k - \lambda }
    \end{equs}
    We note first a convenient representation for $I_n^{(2)}$ in the form of a telescoping sum  
    \begin{equs}
    I_n^{(2)}(\lambda, \lambda_k) &= \frac{1}{\lambda_k - \lambda} \E \left[ \phi_{\lambda, n}(X_{\lambda_k, 1}, \dots, X_{\lambda_k, n}) - \phi_{\lambda, n}(X_{\lambda, 1}, \dots, X_{\lambda, n})\right]  \\ 
    &= \frac{1}{\lambda_k - \lambda} \sum_{j=1}^n \E \left[ \phi_{\lambda, n}(X_{\lambda_k, 1}, \dots, X_{\lambda_k, j}, X_{\lambda, j+1}, \dots,  X_{\lambda, n}) - \phi_{\lambda, n}(X_{\lambda_k, 1}, \dots, X_{\lambda_k, j-1}, X_{\lambda, j}, \dots ,  X_{\lambda, n}) \right] 
    \end{equs}
    We will first argue that we can interchange limits and the sum over $n$. Indeed, we can immediately bound 
    \begin{equs}
    |I_n^{(1)} | &\leq \frac{\norm{f'}_\infty}{u} \\ 
    |I_n^{(2)}| &\leq \frac{1}{\lambda_k - \lambda} \sum_{j=1}^n \E \left[ \norm{\partial_i \phi_{\lambda,n}}_\infty (X_{\lambda_k, j}- X_{\lambda, j})  \right] \leq \frac{n C_{u,v}}{u} \norm{f'}_\infty 
    \end{equs}
    Here we used Lemma \ref{lem:X_lipschitz} to obtain the very last inequality. Hence we see that the map $n \mapsto I_n$ is pointwise bounded by a function which is integrable: 
    \begin{equ}
    \sum_n p_n (1+nC_{u,v}) \frac{\norm{f'}_\infty}{u} =  \left( 1 +C_{u,v} \E Z\right) \frac{\norm{f'}_\infty}{u}  < \infty 
    \end{equ}
    Thus by dominated convergence theorem we have 
    \begin{equ}
    \lim_{k \to \infty} \sum_n p_n I_n(\lambda, \lambda_k ) = \sum_n p_n \lim_{k \to \infty} I_n(\lambda, \lambda_k ) \, , 
    \end{equ}
    granted that the pointwise the limit $\lim_k I_n(\lambda, \lambda_k)$ exists. This we will show now. For $I_n^{(1)}$ this is quite straightforward, since $\mu_{\lambda_k} \xrightharpoonup{} \mu_\lambda$ and $\frac{\phi_{\lambda_k,n} - \phi_{\lambda, n}}{\lambda_{k} - \lambda} \to \partial_\lambda\phi_{\lambda, n}$ uniformly as $f \in C^1([0,1])$. Hence 
    \begin{equ}
    \lim_{k \to \infty} I_n^{(1)}(\lambda, \lambda_k) = \int \partial_\lambda \phi_{\lambda, n} \dd \mu_\lambda^{\otimes n} \, . 
    \end{equ}
    For $I_n^{(2)}$ we extract from the telescoping sum with the fundamental theorem of calculus and the independence of $X_{\lambda_k, 1}, \dots, X_{\lambda, n}$: 
    \begin{equs}
    I^{(2)}_{n}(\lambda, \lambda_k) &= \frac{1}{\lambda_k - \lambda} \sum_{j=1}^n \E \left[ \int_{X_{\lambda, j}}^{X_{\lambda_k, j}} \partial_j \phi_{\lambda, n}(X_{\lambda_k, 1}, \dots, z, \dots, X_{\lambda, n}) \dd z \right] \\ 
    &= -\frac{1}{\lambda} \sum_{j=1}^n \int \partial_j \phi_{\lambda, n} \dd( \mu_{\lambda_k}^{\otimes(j-1)} \otimes R_{\lambda, \lambda_k} \otimes \mu_\lambda^{\otimes(n-j)}) \, . 
    \end{equs}
    Now as $\partial_j \phi_{\lambda, n}$ is itself bounded and continuous due to $f \in C^1([0,1])$, we can pass to the limit, using $ \mu_{\lambda_k} \xrightharpoonup{} \mu_\lambda$ and $R_{\lambda, \lambda_k} \xrightharpoonup{} \tilde{\Gamma}_\lambda$, to obtain 
    \begin{equ}
    \lim_{k \to \infty} I_n^{(2)}(\lambda, \lambda_k) = - \frac{1}{\lambda} \sum_{j=1}^n \int \partial_j \phi_{\lambda, n} \dd (\mu_{\lambda}^{\otimes (j-1)} \otimes \tilde{\Gamma}_{\lambda} \otimes \mu_\lambda^{\otimes(n-j)}) = -\frac{n}{\lambda} \int \partial_1 \phi_{\lambda, n} \dd (\tilde{\Gamma}_{\lambda} \otimes \mu_\lambda^{\otimes (n-1)}) \, . 
    \end{equ}
    The last equality follows from the fact that $\phi_{\lambda, n}$ factors through $y \mapsto \norm{y}_1$ which is permutation invariant. This concludes the proof. 
\end{proof}

\begin{lemma}\label{lem:weak_continuity}
    We have 
    \begin{enumerate}
        \item $\lambda \mapsto \mu_\lambda, (0,m) \to \CM^+$ is weakly continuous. 
        \item $\lambda \mapsto \Gamma_\lambda, (0,m) \to \CM^+$ is weakly continuous. 
    \end{enumerate}
\end{lemma}
\begin{proof}
    \begin{enumerate}
        \item It easy to see that by the Lipschitz continuity of Lemma \ref{lem:X_lipschitz}, and the fact that Lipschitz continuous functions are dense in $C([0,1])$, that we have $\lambda \mapsto \mu_\lambda$ is weakly continuous. 
        \item  Let $\lambda_n \to \lambda$. To prove $\Gamma_{\lambda_n} \rightharpoonup \Gamma_\lambda$ we use the usual subsequence argument. So let $(\lambda_n)_n$ be any subsequence of the original sequence\footnote{using the same indexing by abuse of notation}. First, note that by \eref{eq:Gamma_uniform_bound} and Jensen's inequality, we get a uniform in $\lambda$ bound 
        \begin{equ}
        \norm{\Gamma_\lambda}_{\mathrm{TV}} \leq \lim_{N \to \infty} \norm{\Gamma_\lambda^N}_{\mathrm{TV}} \leq \frac{16}{9} \frac{(\E X_\lambda )  (\E X_\lambda - X_\lambda^2)}{\E X_\lambda^2 } \leq \frac{16}{9} 
        \end{equ}
        Hence, by tightness, we can extract a subsequence such that $(\Gamma_{\lambda_{n_k}})_{k}$ converges weakly in $\CM_\lambda^+$, and denote its limit by $\tilde{\Gamma}_\lambda = \weaklim_{k \to \infty} \Gamma_{\lambda_{n_k}}$. Now for every $k$ it holds 
        \begin{equ}
        \Gamma_{\lambda_{n_k}} = \sigma_{\lambda_{n_k}} + T_{\lambda_{n_k}} \Gamma_{\lambda_{n_k}} \, . 
        \end{equ}
        As clearly $\lambda \mapsto \sigma_\lambda $ is weakly continuous, all that remains is to compute the weak limit of $ T_{\lambda_{n_k}} \Gamma_{\lambda_{n_k}}$. To this end, let $f \in C([0,1])$ and compute $\int f \dd  T_{\lambda_{n_k}} \Gamma_{\lambda_{n_k}}$. We denote by $\Phi_\lambda = \Phi_\lambda^{(m,f)} \in C([0,1]^m;\R)  $ the function such that 
        \begin{equ}
        \int f \dd T_\lambda \Gamma_\lambda = \sum_{m=4}^\infty mp_m \int \Phi_\lambda \dd(\Gamma_\lambda \otimes \mu_\lambda^{\otimes (m-1)}) \, . 
        \end{equ}
        Explicitly, $\Phi_\lambda(x,y_1, \dots, y_{m-1}) = \frac{\lambda}{(\lambda + x + \sum_i y_i )^2}f\left( \frac{x + \sum_i y_i}{\lambda + x + \sum_i y_i} \right) $. Since $\norm{\Phi_{\lambda_{n_k}} - \Phi_\lambda}_\infty \to 0 $ and $\Gamma_{\lambda_{n_k}} \otimes \mu_{\lambda_{n_k}}^{\otimes(m-1)} \xrightharpoonup{} \tilde{\Gamma}_\lambda \otimes \mu_\lambda^{\otimes(m-1)} $ it is easy to see that for every fixed $m \in \N$: 
        \begin{equ}
        \lim_{k \to \infty} \int \Phi_{\lambda_{n_k}} \dd (\Gamma_{\lambda_{n_k}} \otimes \mu_{\lambda_{n_k}}^{\otimes(m-1)}) = \int \Phi_\lambda \dd (\tilde{\Gamma}_\lambda \otimes \mu_\lambda^{\otimes(m-1)}) \, . 
        \end{equ}
        Finally, by inspection 
        \begin{equ}
        \norm{\Phi_\lambda}_\infty \leq \frac{\norm{f}_\infty}{u} \, , 
        \end{equ}
        where $u > 0 $ is some lower bound such that for all $n \in \N:$ $\lambda_n > u$ which we can assume without loss of generality. 
        This, combined with the bound on $\norm{\Gamma_\lambda}_{\mathrm{TV}}$, yields the following uniform in $\lambda_{n_k}$ bound 
        \begin{equ}
        \left| \sum_{m \geq N} mp_m \int  \Phi_{\lambda_{n_k}} \dd (\Gamma_{\lambda_{n_k}} \otimes \mu_{\lambda_{n_k}}^{\otimes(m-1)}) \right| \leq \frac{16}{9} \frac{\norm{f}_\infty}{u} \sum_{m \geq N} mp_m  \, . 
        \end{equ}
        The right-hand-side tends to 0 for $N \to \infty$  as $\E Z = m < \infty$, and hence we can interchange limit taking over $\lambda_{n_k} \to \lambda$ and summing over $m \in \N$. Therefore $T_{\lambda_{n_k}} \Gamma_{\lambda_{n_k}} \xrightharpoonup{} T_\lambda \tilde{\Gamma}_\lambda $, and thus 
        \begin{equ}
        \tilde{\Gamma}_\lambda = \sigma_\lambda + T_\lambda \tilde{\Gamma}_\lambda, 
        \end{equ}
        from which, by $\tilde{\Gamma}_\lambda \in \CM^+_\lambda$ and uniqueness follows $\tilde{\Gamma}_\lambda = \Gamma_\lambda$. So by the subsequence argument, we have proven that $\Gamma_{\lambda_n} \xrightharpoonup{} \Gamma_{\lambda}$. 
    \end{enumerate}
\end{proof}

\newpage 

\begin{lemma}\label{lem:omega_positive_measure}
    The measure $\Omega_\lambda$ is positive. 
\end{lemma}
\begin{proof}
    Since $\Omega_\lambda = \sigma_\lambda - (1 - \lambda) \Gamma_\lambda$, and $\sigma_\lambda, \Gamma_\lambda \geq 0$, positivity is clear when $\lambda \in [1,m)$. So let $\lambda \in (0,1)$. Then observe that 
    \begin{equs}
    \Omega_\lambda = \sigma - (1-\lambda)(\sigma + T \Gamma) = \lambda \sigma + \lambda T\Gamma - T \Gamma  = \lambda \Gamma - T \Gamma = \sum_{j=0}^\infty T^{j+1}(\lambda \alpha - \sigma) \, . 
    \end{equs} 
    Now for $\lambda <1$ one sees that with $(\lambda \alpha - \sigma)(\dd x) = x (\lambda  - (1-x)) \mu(\dd x) $, that as $\supp \mu_\lambda \subset [(1-\lambda / 4)_+, 1]$, the elementary inequality 
    \begin{equ}
    x+ (\lambda - 1) \geq 1-\lambda / 4 + (\lambda - 1) = \frac{3}{4} \lambda > 0 
    \end{equ}
    yields $\lambda \alpha - \sigma \geq 0$. As $T$ is a positive operator, positivity of $\Omega_\lambda$ follows. 
\end{proof}

\paragraph{Acknowledgements} FS is funded by the European Research Council (ERC) under
the European Union's Horizon 2020 research and innovation programme (Grant agreement No.
101045082), funded also by the Deutsche Forschungsgemeinschaft (DFG, German Research Foundation) under Germany's Excellence Strategy EXC 2044/2 –390685587, Mathematics Münster: Dynamics–Geometry–Structure, and has been supported by the Research Training Group “Rigorous Analysis of
Complex Random Systems” (RTG 3027) funded by the DFG (Project
Number 524444762). The author is grateful to Jan Nagel for helpful discussions and for first bringing this theorem to their attention.

\paragraph{Tool usage statement.} OpenAI's ChatGPT 6 Astra was used during the initial stages of this work. The proof in its final form was developed by the author, who takes full responsibility for all claims and arguments made in the paper. 

\bibliographystyle{myalpha} 
\bibliography{refs} 


\small 
\vspace{2em}
\noindent
(F.~Schwarzfischer) \textsc{Institut für Analysis und Numerik, University of Münster, Münster (Westf),  Germany}\\
\hspace*{1.5em}\textit{Email address}: \texttt{felix.schwarzfischer@uni-muenster.de}

\end{document}